\documentclass[a4paper,10pt]{article}
\usepackage{paper-en}
\usepackage{hyperref}

\def\thetitle{Veronese minimizes normal curvatures}
\def\theauthors{Anton Petrunin}

\hypersetup{colorlinks=true,
citecolor=black,
linkcolor=black,
anchorcolor=black,
filecolor=black,
menucolor=black,
urlcolor=black,
pdftitle={\thetitle},
pdfauthor={\theauthors}
}

\begin{document}

\title{\thetitle}
\author{\theauthors
\blfootnote{This work was partially supported by the National Science Foundation, grant DMS-2005279.}}
\date{}
\maketitle

\begin{abstract}
Suppose \( M \) is a closed submanifold in a Euclidean ball of sufficiently large dimension.
We give an optimal bound on the normal curvatures, guaranteeing that \( M \) is a sphere.
The border cases consist of Veronese embeddings of the four projective planes.
\end{abstract}

\paragraph{Introduction.}
Let \( M \subset \mathbb{R}^d \) be a closed smooth \( n \)-dimensional submanifold.
Assume that \( d \) is large and \( M \) is contained in a ball of radius \( r \).
\textit{What can be said about the maximal normal curvature of \( M \)?}
In other words, \textit{what can be said about the maximal curvature of geodesics in \( M \) when considered as curves in \( \mathbb{R}^d \)?}

First, note that the curvatures cannot be smaller than \( \tfrac{1}{r} \) at all points.
Moreover,
\textit{the average value of \( |H| \) must be at least \( n \cdot \tfrac{1}{r} \)} \cite[28.2.5]{burago-zalgaller}, \cite[3.1]{petrunin2024a};
here \( H \) denotes the mean curvature vector.
This statement is a straightforward generalization of István Fáry's result about closed curves in a ball \cite{fary,tabachnikov}.

On the other hand, the \(n\)-dimensional torus can be embedded into an \(r\)-ball with all normal curvatures equal to \( \sqrt{3 \cdot n / (n + 2)} \cdot \tfrac{1}{r} \).
This embedding was discovered by Michael Gromov \cite[2.A]{gromov3}, \cite[1.1.A]{gromov2}.
This bound is optimal; that is, \textit{any smooth \(n\)-dimensional torus in an \(r\)-ball must have normal curvature at least \( \sqrt{3 \cdot n / (n + 2)} \cdot \tfrac{1}{r} \) at some point}~\cite{petrunin2024a}.

Gromov's examples easily imply the following:
\textit{any closed smooth manifold \( M \) admits a smooth embedding into an \(r\)-ball of sufficiently large dimension with normal curvatures less than \( \sqrt{3} \cdot \tfrac{1}{r} \)}
\cite[1.D]{gromov3}, \cite[1.1.C]{gromov2}.

\textit{But what happens between \( \tfrac{1}{r} \) and \( \sqrt{3} \cdot \tfrac{1}{r} \)?}

In this note, we consider embeddings into an \(r\)-ball with normal curvatures at most \( \tfrac{2}{\sqrt{3}} \cdot \tfrac{1}{r} \).
We show that if the inequality is strict, the manifold must be homeomorphic to a sphere (see §~\ref{thm:strict}).
For the non-strict inequality, in addition to spheres, the possible cases include the real, complex, quaternionic, and octonionic planes mapped by rescaled Veronese embeddings (see §~\ref{thm:=}).

\paragraph{Sphere theorem.}
\label{thm:strict}
\textit{Let \( M \) be a closed smooth \( n \)-dimensional submanifold contained in a closed \( r \)-ball in \( \mathbb{R}^d \).
Suppose the normal curvatures of \( M \) are strictly less than \( \tfrac{2}{\sqrt{3}} \cdot \tfrac{1}{r} \).
Then \( M \) is homeomorphic to the \( n \)-sphere.}

\medskip

In the proof we will use the following version of Axel Schur's lemma \cite{shur}.

\begin{thm}{Bow lemma}
Let $\gamma_1\:[a,b]\to\mathbb{R}^2$ and $\gamma_2\:[a,b] \to\mathbb{R}^d$ be two $C^{1,1}$-smooth unit-speed curves.
Suppose that curvature of $\gamma_1$ does not exceed the curvature of $\gamma_2$ at any time moment
and the curve
$\gamma_1$ is an arc of a convex curve; that is, it runs in the boundary of a convex plane figure.
Then the distance between the endpoints of $\gamma_1$ cannot exceed the  distance between the endpoints of $\gamma_2$; that is,
\[|\gamma_1(b)-\gamma_1(a)|\leqslant |\gamma_2(b)-\gamma_2(a)|.\]

Moreover, in case of equality, $\gamma_2=\iota \circ \gamma_1$ for some isometry $\iota$ from $\mathbb{R}^2$ to a plane in $\mathbb{R}^d$.
\end{thm}

The first part of the lemma follows from \cite[Theorem 5.1]{sullivan}, and the equality case follows from its proof.

\parit{Proof of the sphere theorem.}
We may assume that \( n \geqslant 2 \);
otherwise, there is nothing to prove.

Denote the \( r \)-ball by \( \mathbb{B}^d \).
We can assume \( r = \tfrac{1}{\sqrt{3}} \), meaning \( r \) is the circumradius of an equilateral triangle with unit side.
Consequently, the normal curvatures of \( M \) are less than 2.

Choose a unit-speed geodesic \( \gamma\colon  [0, \tfrac{\pi}{2}] \to M \);
let \( x = \gamma(0) \) and \( y = \gamma(\tfrac{\pi}{2}) \).
By assumption, the curvature of \( \gamma \) in \( \mathbb{R}^d \) is less than 2.
Applying the bow lemma, we obtain \( |x - y| > 1 \).

Let $\Pi$ be the perpendicular bisector of $[x,y]$.
Since the curvature of $\gamma$ is less than $2$,
\[
\measuredangle(\gamma'(t_0), \gamma'(t)) < 2 \cdot |t - t_0|
\quad\text{and}\quad
\langle \gamma'(t_0), \gamma'(t) \rangle > \cos(2 \cdot |t - t_0|),
\]
if $t \ne t_0$.
Therefore,
\[
\langle y - x, \gamma'(t_0) \rangle > \int\limits_0^{\frac{\pi}{2}} \cos(2 \cdot |t - t_0|) \cdot dt \geqslant 0.
\]
In particular, the derivative of the function
$f \colon t \mapsto \langle y - x, \gamma(t) \rangle$
is positive.
Therefore, $\gamma$ intersects $\Pi$ transversely at a single point; denote this point by $s$.

Choose a unit vector $\mathsc{u} \in \T_x$.
Let $\gamma_{\mathsc{u}} \colon [0, \tfrac{\pi}{2}] \to M$ be the unit-speed geodesic starting from $x$ in the direction of $\mathsc{u}$, and let $z = \gamma_{\mathsc{u}}(\tfrac{\pi}{2})$.
The argument above shows that $|x - z| > 1$.

\begin{wrapfigure}{r}{36 mm}
\vskip-4mm
\centering
\includegraphics{mppics/pic-10}
\vskip0mm
\end{wrapfigure}

Denote by $H_x$ and $H_y$ the closed half-spaces bounded by $\Pi$ containing $x$ and $y$, respectively.
Assume $z\in H_x$.
Then
$|y - z| \geqslant |x - z| > 1$.
Since $|x - y| > 1$, the triangle $[xyz]$ has all sides greater than $1$, which is impossible because $x, y, z \in \mathbb{B}^d$.
Therefore, $\gamma_{\mathsc{u}}$ meets $\Pi$ before $\tfrac{\pi}{2}$; denote by $r(\mathsc{u})$ the first such time.

Let us show that the function $\mathsc{u} \mapsto r(\mathsc{u})$ is smooth.
In other words, $\gamma_{\mathsc{u}}$ intersects $\Pi$ transversely at time $r(\mathsc{u})$.
Assume this is not the case, then $\gamma_{\mathsc{u}}$ is tangent to $\Pi$ at $r(\mathsc{u})$.
Let $\hat{\gamma}_{\mathsc{u}}$ be the concatenation of the reflection of $\gamma_{\mathsc{u}}|_{[0, r(\mathsc{u})]}$ across $\Pi$ and $\gamma_{\mathsc{u}}|_{[r(\mathsc{u}), \frac{\pi}{2}]}$.
Note that $\hat{\gamma}_{\mathsc{u}}$ is $C^1$-smooth and $C^\infty$-smooth everywhere except at $r(\mathsc{u})$.
Therefore, the bow lemma applies to $\hat{\gamma}_{\mathsc{u}}$, and hence $|y - z| \z> 1$.
Again, all sides of the triangle $[xyz]$ are greater than $1$; hence, it cannot lie in $\mathbb{B}^d$ --- a contradiction.

It follows that the set
\[
V_x = \set{t \cdot \mathsc{u} \in \T_x}{|\mathsc{u}| = 1, \quad 0 \leqslant t \leqslant r(\mathsc{u})}
\]
is diffeomorphic to the closed \(n\)-disc.
Denote by \(W_x\) the connected component of \(M \cap H_x\) that contains \(x\).

By the Gauss formula \cite[Lemma 5]{petrunin2024}, the sectional curvatures of \(M\) are less than \(4\).
In particular, the exponential map \(\exp_x \colon \T_x \to M\) is a local diffeomorphism within the \(\tfrac{\pi}{2}\)-ball centered at the origin in \(\T_x\).
It follows that the exponential map defines a local diffeomorphism $V_x \to W_x$.
In particular, $W_x$ is a smooth manifold with boundary.
Since $V_x$ is simply connected, $\exp_x$ defines a diffeomorphism between $V_x$ and $W_x$.
Thus, $W_x$ is a closed topological $n$-disc, and $\partial W_x$ is a smooth hypersurface in $M$.

Let us swap the roles of \(x\) and \(y\) and repeat the construction.
We obtain another closed topological \(n\)-disc \(W_y \subset M\), bounded by a smooth hypersurface \(\partial W_y\).

Observe that $\partial W_x$ intersects $\partial W_y$ at $s$.
Furthermore, both $\partial W_x$ and $\partial W_y$ are connected components of $s$ in $M\cap \Pi$.
Therefore, $\partial W_x=\partial W_y$.
That is, $M$ can be obtained by gluing two $n$-discs by a diffeomorphism between their boundaries.
Hence $M$ is homeomorphic to the $n$-sphere.
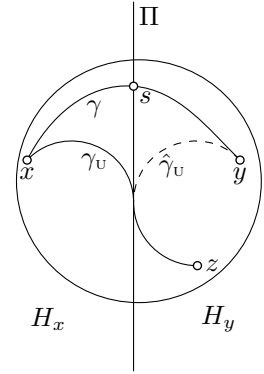
\qeds

\paragraph{Veronese embeddings.}\label{veronese}
The real, complex, and quaternionic projective spaces of dimension \(n\), along with the octonionic projective plane, are denoted by \(\RP^n\), \(\CP^n\), \(\HP^n\), and \(\OP^2\), respectively.
We assume that each of these spaces is equipped with the canonical metric; in particular, all the spaces have closed geodesics of length~\(\pi\).

\begin{thm}{Proposition}
There are smooth isometric embeddings
\begin{itemize}
 \item $\RP^n \hookrightarrow\mathbb{R}^d$ for $d\geqslant n+\tfrac12\cdot n\cdot(n+1)$;
 \item $\CP^n \hookrightarrow\mathbb{R}^d$ for $d\geqslant n+n\cdot(n+1)$;
 \item $\HP^n \hookrightarrow\mathbb{R}^d$ for $d\geqslant n+2\cdot n\cdot(n+1)$;
 \item $\OP^2 \hookrightarrow\mathbb{R}^d$ for $d\geqslant 26$;
\end{itemize}
that map each geodesic to a round circle.

Moreover,

\begin{subthm}{2}
all normal curvatures of the images under these embeddings are equal to~2,
\end{subthm}

\begin{subthm}{r}
the images of these embeddings lie on a sphere of radius \( r_n = \sqrt{n / (2 \cdot n + 2)} \)
(for \( \OP^2 \), we assume \( n = 2 \)).
\end{subthm}

\end{thm}

This proposition, along with a detailed proof, is given in Kunio Sakamoto's paper \cite[§ 2]{sakamoto}.
These embeddings will be referred to as \emph{Veronese embeddings}.
Note that \( r_n \) is the circumradius of a regular \( n \)-simplex with unit edge length, and \(r_2 = 1 / \sqrt{3}\).
Thus, the second part of the proposition implies that our sphere theorem is optimal.

The Veronese embeddings have an explicit algebraic description and possess several remarkable geometric properties.
In particular, these embeddings are equivariant, and their images are minimal submanifolds in the \( r_n \)-spheres.
All of these properties are discussed in the cited paper by Kunio Sakamoto.

We will need the following claim, which is a special case of John Little's theorem \cite{little} and the more general results of Kunio Sakamoto \cite[Theorem 3]{sakamoto}.

\begin{thm}{Claim}
Let \( M \) be a closed  submanifold in \( \mathbb{R}^d \).
Suppose that every geodesic in \( M \) is a planar circle of radius \( \tfrac{1}{2} \).
Then \( M \) is either a round sphere of radius \( \tfrac{1}{2} \) or one of the Veronese embeddings described above.
\end{thm}

\paragraph{Rigidity theorem.}\label{thm:=}
\textit{Let \( M \) be a closed smooth \( n \)-dimensional submanifold in a closed \( r \)-ball in \( \mathbb{R}^d \).
Suppose the normal curvatures of \( M \) are at most \( \tfrac{2}{\sqrt{3}} \cdot \tfrac{1}{r} \).
If \( M \) is not homeomorphic to a sphere, then up to rescaling, it is congruent to the image of a Veronese embedding of a projective plane: \( \RP^2 \), \( \CP^2 \), \( \HP^2 \), or~\( \OP^2 \).}

\medskip

\parit{Proof.}
Assume that \( M \) is not homeomorphic to a sphere;
in this case, \( n \geqslant 2 \).
As before, let \( \mathbb{B}^d \) denotes the \( r \)-ball in \( \mathbb{R}^d \), where we assume \( r = \tfrac{1}{\sqrt{3}} \);
therefore, the normal curvatures of \( M \) are at most 2.

By the proposition in §~\ref{veronese}, the images of the Veronese embedding of the projective planes satisfy the assumptions of the theorem.
It remains to show that no other embeddings of this type exist.

Choose a unit-speed geodesic \( \gamma \colon [0, \tfrac{\pi}{2}] \to M \);
let \( x = \gamma(0) \) and \( y = \gamma(\tfrac{\pi}{2}) \).
The argument in our sphere theorem implies that \( |x-y| = 1 \).
The rigidity case in the bow lemma implies that \( \gamma \) is a half-circle of curvature 2.
Since any two points in \( M \) can be connected by a geodesic, we conclude that all geodesics in \( M \) are circles of curvature 2 in \( \mathbb{R}^d \).

Since \( M \) is not a sphere, the claim in §~\ref{veronese} implies that \( M \) is given by one of the Veronese embeddings of \( \RP^n \), \( \CP^n \), \( \HP^n \) for \( n \geqslant 2 \), or \( \OP^2 \).

It remains to show that \( n \leqslant 2 \).
Assume the contrary.
Note that if $n\geqslant 3$, then each space \( \RP^n \), \( \CP^n \), and \( \HP^n \) contains four points at a distance \( \tfrac{\pi}{2} \) from one another.
The corresponding points would lie at a distance 1 in \( \mathbb{B}^d \), which is impossible.
\qeds

\paragraph{Final remarks.}
Recall that the Veronese embeddings map \( \RP^n \), \( \CP^n \), and \( \HP^n \) into balls of radius \( r_n = \sqrt{n/(2\cdot n + 2)} \), which is the circumradius of a regular \( n \)-simplex with unit edge length.

This note is motivated by the following question; see \cite[1.4]{petrunin2024a} and \cite{petrunin2023}.

\begin{thm}{Question}
Is it true that the Veronese embedding minimizes the maximum normal curvature among all smooth embeddings of \(\RP^n\) into an \(r_n\)-ball in a Euclidean space of sufficiently large dimension?
\end{thm}

The same question can be posed for \( \CP^n \) and \( \HP^n \).
A keen reader might have noticed that the case $n=2$ is already solved.

Gromov's construction provides numerous smooth immersions of the torus \( \TT^n \) into \( r \)-balls with all normal curvatures \( \sqrt{3\cdot n/(n+2)}\cdot\tfrac{1}{r}\).
The proof of the theorem in \cite{petrunin2024a} implies that such immersions map the torus to the boundary of the ball;
it also precisely describes its second fundamental form.
Nevertheless, the following question remains open.

\begin{thm}{Question}
Can one characterize all smooth immersions of the $n$-torus into $r$-balls
with all normal curvatures equal to $\sqrt{3\cdot n/(n+2)}\cdot\tfrac{1}{r}$?
\end{thm}

\begin{thm}{Question}
Let \( M \) be as in our sphere theorem.
Must it be diffeomorphic to the standard \( n \)-sphere?
\end{thm}

I suspect that the answer to the previous question is \textit{yes}.
If $M$ lies on the boundary of the $r$-ball, then by the Gauss formula \cite[Lemma 5]{petrunin2024}, $M$ has strictly quarter-pinched curvature.
Therefore, in this case, the answer is \textit{yes} by the differentiable sphere theorem \cite{brendle-schoen}.

\paragraph{Acknowledgments.}
I would like to thank Alexander Lytchak and Ricardo Mendes for their help.
The proof of the rigidity theorem was suggested by Ricardo Mendes.
My original proof was based on the diameter-rigidity theorem \cite{gromoll-grove,wilking}.

This work was partially supported by the National Science Foundation, grant DMS-2005279.

{\sloppy
\def\emph{\textit}
\printbibliography[heading=bibintoc]

@book {burago-zalgaller,
    AUTHOR = {Burago, Yu. D. and Zalgaller, V. A.},
     TITLE = {Geometric inequalities},
    SERIES = {Grundlehren der mathematischen Wissenschaften},
    VOLUME = {285},
% PUBLISHER = {Springer-Verlag, Berlin},
      YEAR = {1988},
   %  PAGES = {xiv+331},
      ISBN = {3-540-13615-0},
   MRCLASS = {52A40 (53-02)},
  MRNUMBER = {936419},
       DOI = {10.1007/978-3-662-07441-1},
       URL = {https://doi.org/10.1007/978-3-662-07441-1},
}

@Article{fary,
    Author = {Fáry, I.},
    Title = {Sur certaines inégalités géométriques},
    FJournal = {Acta Scientiarum Mathematicarum},
    Journal = {Acta Sci. Math.},
    ISSN = {0001-6969},
    Volume = {12},
    Pages = {117--124},
    Year = {1950},
 %   Publisher = {University of Szeged, Bolyai Institute, Szeged},
    Language = {French},
    Zbl = {0039.16801}
}

@InCollection{tabachnikov,
    Author = {S. {Tabachnikov}},
    Title = {The tale of a geometric inequality.},
    BookTitle = {{MASS selecta: teaching and learning advanced undergraduate mathematics}},
    ISBN = {0-8218-3363-4/hbk},
    Pages = {257--262},
    Year = {2003},
   % Publisher = {Providence, RI: American Mathematical Society (AMS)},
    Language = {English},
    MSC2010 = {53A04 53C65},
    Zbl = {1093.53001}
}

@article{gromov2,
doi = {10.48550/ARXIV.2210.13256},
author = {Gromov, M.},
title = {Curvature, Kolmogorov diameter, Hilbert rational designs and overtwisted immersions},
year={2022},
      eprint={2210.13256},
      archivePrefix={arXiv},
      primaryClass={math.DG}
}

@misc{gromov3,
author = {Gromov, M.},
title = {Isometric immersions with controlled curvatures},
year={2022},
      eprint={2212.06122v1},
      archivePrefix={arXiv},
      primaryClass={math.DG}
}

@article {little,
    AUTHOR = {Little, J. A.},
     TITLE = {Manifolds with planar geodesics},
   JOURNAL = {J. Differential Geometry},
  FJOURNAL = {Journal of Differential Geometry},
    VOLUME = {11},
      YEAR = {1976},
    NUMBER = {2},
     PAGES = {265--285},
      ISSN = {0022-040X,1945-743X},
   MRCLASS = {53C40},
  MRNUMBER = {417992},
MRREVIEWER = {O.\ Kowalski},
       URL = {http://projecteuclid.org/euclid.jdg/1214433424},
}

@article {shur,
    AUTHOR = {Schur, A.},
     TITLE = {\"{U}ber die {S}chwarzsche {E}xtremaleigenschaft des {K}reises
              unter den {K}urven konstanter {K}r\"{u}mmung},
   JOURNAL = {Math. Ann.},
  FJOURNAL = {Mathematische Annalen},
    VOLUME = {83},
      YEAR = {1921},
    NUMBER = {1-2},
     PAGES = {143--148},
      ISSN = {0025-5831},
   MRCLASS = {DML},
  MRNUMBER = {1512005},
       DOI = {10.1007/BF01464234},
       URL = {https://doi.org/10.1007/BF01464234},
}

@article {gromoll-grove,
    AUTHOR = {Gromoll, D. and Grove, K.},
     TITLE = {A generalization of {B}erger's rigidity theorem for positively
              curved manifolds},
   JOURNAL = {Ann. Sci. \'Ecole Norm. Sup. (4)},
  FJOURNAL = {Annales Scientifiques de l'\'Ecole Normale Sup\'erieure.
              Quatri\`eme S\'erie},
    VOLUME = {20},
      YEAR = {1987},
    NUMBER = {2},
     PAGES = {227--239},
      ISSN = {0012-9593},
   MRCLASS = {53C20},
  MRNUMBER = {911756},
MRREVIEWER = {Viktor\ Schroeder},
       URL = {http://www.numdam.org/item?id=ASENS_1987_4_20_2_227_0},
}

@article {wilking,
    AUTHOR = {Wilking, B.},
     TITLE = {Index parity of closed geodesics and rigidity of {H}opf
              fibrations},
   JOURNAL = {Invent. Math.},
  FJOURNAL = {Inventiones Mathematicae},
    VOLUME = {144},
      YEAR = {2001},
    NUMBER = {2},
     PAGES = {281--295},
      ISSN = {0020-9910,1432-1297},
   MRCLASS = {53C12 (53C22 53C24)},
  MRNUMBER = {1826371},
       DOI = {10.1007/PL00005801},
       URL = {https://doi.org/10.1007/PL00005801},
}

@article{petrunin2024,
      title={Tubed embeddings}, 
      author={A. Petrunin},
      year={2024},
      eprint={2402.16195},
      archivePrefix={arXiv},
      primaryClass={math.DG},
      url={https://arxiv.org/abs/2402.16195}, 
}

@MISC {petrunin2023,
    TITLE = {Normal curvature of Veronese embedding},
    AUTHOR = {A. Petrunin},
    HOWPUBLISHED = {MathOverflow},
   % NOTE = {URL:https://mathoverflow.net/q/445819 (version: 2024-06-10)},
    EPRINT = {https://mathoverflow.net/q/445819},
    URL = {https://mathoverflow.net/q/445819}
}

@article {sakamoto,
    AUTHOR = {Sakamoto, K.},
     TITLE = {Planar geodesic immersions},
   JOURNAL = {Tohoku Math. J. (2)},
  FJOURNAL = {The Tohoku Mathematical Journal. Second Series},
    VOLUME = {29},
      YEAR = {1977},
    NUMBER = {1},
     PAGES = {25--56},
      ISSN = {0040-8735,2186-585X},
   MRCLASS = {53C40},
  MRNUMBER = {470913},
MRREVIEWER = {O.\ Kowalski},
       DOI = {10.2748/tmj/1178240693},
       URL = {https://doi.org/10.2748/tmj/1178240693},
}

@article {petrunin2024a,
    AUTHOR = {Petrunin, A.},
     TITLE = {Gromov's tori are optimal},
   JOURNAL = {Geom. Funct. Anal.},
  FJOURNAL = {Geometric and Functional Analysis},
    VOLUME = {34},
      YEAR = {2024},
    NUMBER = {1},
     PAGES = {202--208},
      ISSN = {1016-443X,1420-8970},
   MRCLASS = {53C42},
  MRNUMBER = {4706446},
       DOI = {10.1007/s00039-024-00663-0},
       URL = {https://doi.org/10.1007/s00039-024-00663-0},
}

@article {brendle-schoen,
    AUTHOR = {Brendle, S. and Schoen, R.},
     TITLE = {Manifolds with 1/4-pinched curvature are space forms},
   JOURNAL = {J. Amer. Math. Soc.},
  FJOURNAL = {Journal of the American Mathematical Society},
    VOLUME = {22},
      YEAR = {2009},
    NUMBER = {1},
     PAGES = {287--307},
      ISSN = {0894-0347,1088-6834},
   MRCLASS = {53C20 (53C44)},
  MRNUMBER = {2449060},
MRREVIEWER = {Fr\'ed\'eric\ Robert},
       DOI = {10.1090/S0894-0347-08-00613-9},
       URL = {https://doi.org/10.1090/S0894-0347-08-00613-9},
}

@incollection {sullivan,
    AUTHOR = {Sullivan, J.},
     TITLE = {Curves of finite total curvature},
 BOOKTITLE = {Discrete differential geometry},
    SERIES = {Oberwolfach Semin.},
    VOLUME = {38},
     PAGES = {137--161},
 PUBLISHER = {Birkh\"auser, Basel},
      YEAR = {2008},
      ISBN = {978-3-7643-8620-7},
   MRCLASS = {53A04},
  MRNUMBER = {2405664},
MRREVIEWER = {Jesse\ Ratzkin},
       DOI = {10.1007/978-3-7643-8621-4\_7},
       URL = {https://doi.org/10.1007/978-3-7643-8621-4_7},
}
\fussy
}
\end{document}